# Piecewise Linearization of Quadratic Branch Flow Limits by Irregular Polygon

Parikshit Pareek and Ashu Verma, *Senior Member, IEEE*

*Abstract*— **This letter addresses the issue of linearization of quadratic thermal limits of transmission lines for linear OPF formulation. A new irregular polygon based linearization is proposed. The approach is purely based on geometrical concepts and does not introduce any optimization problem. Comparison of the number of constraints is given with different errors for different branch capabilities and systems. Test case analysis is also presented to validate the irregular polygon linearization strategy with optimal value and computational time results.**

*Index Terms*—Piecewise Linearization, Linear Optimal Power Flow (LOPF).

## I. Piecewise Linearization

The thermal limit linearization using the inner approximation (conservative approach) of the circle is presented in reference [1] for Optimal Power Flow (OPF). It leaves a small area of circle unutilized as the error as depicted in Fig. (1). The $e_{max}$ is the maximum perpendicular distance between the line and circumference of the circle occurring at the mid-point $(-Q_{ij}, P_{ij})$ of that line, called Sagitta.

For a $n^{th}$ line segment of length $L_{n,n+1}$, for branch MVA limit $S_i$, the Sagitta (referred as error hereafter) is:

$$e_{max,n} = e_{m,n} = S_i - \sqrt{S_i^2 - (L_{n,n+1}/2)^2} \quad (1)$$

A relation between $L_{n,n+1}$ and the angle difference $(\Delta\theta_n)$ of $n^{th}$ and $(n+1)^{th}$ point:

$$L_{n,n+1} = \sqrt{2}S_i\sqrt{1 - \cos(\Delta\theta_n)} \quad (2)$$

By replacing $L_{n,n+1}$ in the relation of $e_{m,n}$:

$$\Delta\theta_n = \cos^{-1}\left\{2(1 - e_{m,n}/S_i)^2 - 1\right\} \quad (3)$$

For regular polygon based linearization proposed in the literature:

$$L_{n,n+1} = L_{reg}; \; e_{m,n} = e_{m,reg} \quad n = 1,2,\ldots,M_{reg} \quad (4)$$

$$\Delta\theta_n = \Delta\theta_{reg} = 2\pi/M_{reg} \quad n = 1,2,\ldots,M_{reg} \quad (5)$$

$$M_{reg} = 2\pi/\cos^{-1}\left\{2(1 - e_{m,n}/S_i)^2 - 1\right\} \quad (6)$$

The value of $\theta$, ranging 0 to $2\pi$, is divided in $M_{reg}$ same length line segments. Some strategies are proposed in the literature for the reduction of circle area, that excludes the area inscribed by the arc having angle $\alpha$ in each of the quadrant [2], [3].

Previous work does not reflect upon the $M_{reg}$ and $\alpha$ selection and takes arbitrary value by experience or trials. The $\alpha$, as used in [2], [3] and shown in Fig. 1, is the angle by which the area of circle is reduced in each quadrant.

### A. Proposed Irregular Polygon Linearization

As heavily loaded network branches operate with high $(P_{ij}/Q_{ij})$ ratio, the minimum error is desired to be there. While, in low $P_{ij}$ regions, the branches will be far from limit MVA, i.e., linearization accuracy can be low. This gives the possibility of accurate linearization with a small number of constraints. Mathematically,

$$e_{m,n} = f(\Delta\theta_n) \quad (7)$$

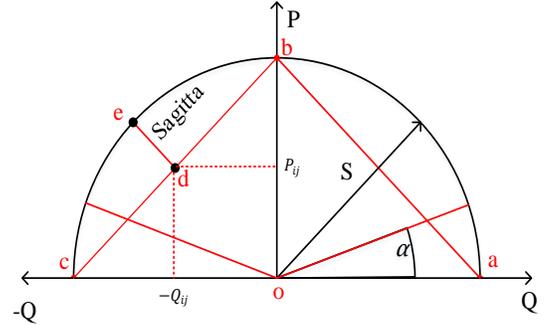

Fig. 1. Sagitta and $\alpha$ in P-Q plane for piecewise linearization of the circle

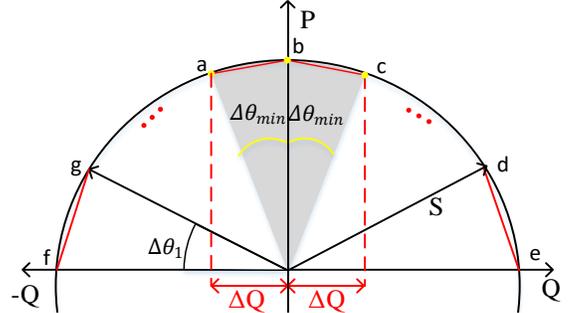

Fig. 2 Depiction of selection criteria for $\Delta\theta_{min}$ and $\Delta Q$

Where, $\quad \Delta\theta_n \neq constant \quad (8)$

Therefore, an irregular polygon with $M_{irr}$ sides is made which linearizes the circle with variable errors. The objectives of using irregular polygon are: 1) minimum error near P-axis and 2) higher line length ($L$) near the Q-axis. To achieve these, the property of circle having variable slope is used and Q-axis is now divided into equal segments i.e. $\Delta Q$. The minimum error is in the line segments on the both sides of P-axis where we want more branches to be operated. This in turn allows a high-power factor branch to utilize its full capacity. Hence, the minimum value of $\Delta\theta_n$ should be calculated accordingly as in Fig. (2). If maximum error desired in the line segment from '$a$' to '$b$' and '$b$' to '$c$' is $e_{m,min}$ then from (3):

$$\Delta\theta_{min} = \cos^{-1}\left\{2\left(1 - e_{m,min}/S_i\right)^2 - 1\right\} \quad (9)$$

$$\Delta Q = S_i * \sin(\Delta\theta_{min}) \quad (10)$$

$$mq = S_i/\Delta Q = \sin^{-1}(\Delta\theta_{min}) \quad (11)$$

Let,
$$\Delta\theta_{min} = \cos^{-1}(\varphi) \quad (12)$$
$$M_{irr} = 4/\sin(\cos^{-1}\{\varphi\}) \quad (13)$$

Using right-angle triangle properties,
$$\cos^{-1}\{\varphi\} = \sin^{-1}\left\{\sqrt{1-\varphi^2}\right\} \quad (14)$$
$$M_{irr} = 4*mq = 4/\sqrt{1-\varphi^2} \quad (15)$$

Replacing $\varphi$ and expanding the quadratic term:
$$M_{irr} = 2/(1-e_{m,min}/S_i)\sqrt{1-(1-e_{m,min}/S_i)^2} \quad (16)$$

Here, $mq$ is the number of polygon sides in a quadrant, $i$ is branch index with $e_{m,min} = \min\{e_{m,n}\}$ for $n = 1, 2 \dots M_{irr}$. As $mq$ is rounded off to next integer, the actual value of $e_{m,min}$ is always less then $e_{m,min}$ used in (16). Both conceptually and figuratively it is clear that the value of $\Delta\theta_1$ can influence the solution of OPF problem (fig. 2). Therefore, it is imperative to check $\Delta\theta_1$, $L_{fg}$ and $e_{m,1}$:

$$L_{fg} = \sqrt{(\Delta Q)^2 + (\Delta P)^2} \quad (17)$$

Here, $\Delta P = P_g = (S_i^2 - Q_g^2)^{1/2}$; $Q_g = -S_i + \Delta Q$ from (11):
$$L_{fg} = \sqrt{2S_i\Delta Q} = S_i\sqrt{2/mq} \quad (18)$$

From (2) and (19) for a branch having '$S_i$' as MVA limit:
$$\Delta\theta_1 = \cos^{-1}(1 - 1/mq) \quad (19)$$

The effect of higher values of $\Delta\theta_1$ on the reactive power flow is not as much as of $\alpha$ in [2]. Nevertheless, it is essential to keep the value of $\Delta\theta_1$ low to ensure OPF solution accuracy and convergence. A relation between $e_{m,n}$ and $\Delta\theta_n$ is obtained from (1) and (2) as:

$$e_{m,n} = 2S_i * \sin^2(\Delta\theta_n/4) \quad (20)$$

Linearization can also be done by selecting a fraction of error as minimum value desired, but it should be avoided as fraction or percentage error does not reflect the actual situation of error uniformly (1% for 1000=10 and for 16=0.16).

The Fig. (3a) shows the error associated with different sides of an irregular polygon with $S_i$ as 16 MVA and $M_{irr} = 40$. The side number is starting from positive Q-axis and increasing in an anti-clockwise direction. The minimum error is with side number 10,11 and 30,31 which are adjoining to positive and negative P-axis respectively. Further, Fig. (3b) shows side length variations for the same irregular polygon. The maximum length is at 1,20,21, and 40 number side which are adjacent to the Q-axis. Thus, the proposed approach has been able to meet the two objectives, minimum error near P-axis and maximum side length near Q-axis.

### B. General Strategy for Hot-start Algorithms

For hot-start or iterative OPF procedures, the linearization can be done near the current operating status. The letter proposes a strategy to introduce the irregular polygon linearization for non-iterative Linear Optimal Power Flow formulations where the operating status of the system is not known prior to the solution. The strategy to linearize at an angle $\theta$, for desired approximation error, is given in the Fig.4.

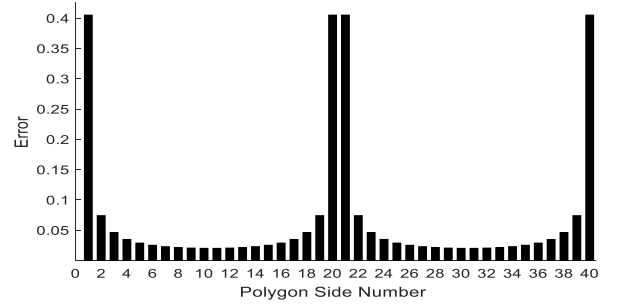

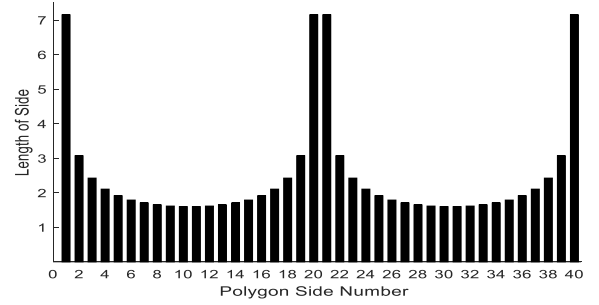

Fig. 3. (a) Error ($e_{m,n}$) and (b) length of irregular polygon sides

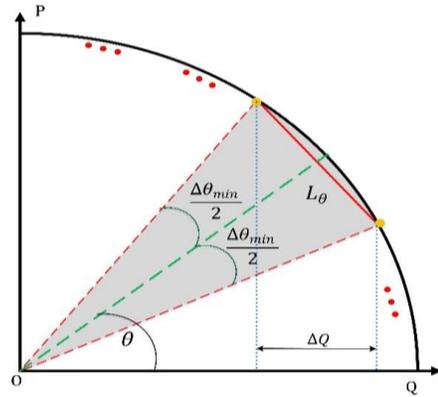

Fig. 4. General strategy for linearization at any angle $\theta$

Let, the angle $\theta$ is the midpoint of a line segment $L_\theta$ and inscribing the angle $\Delta\theta_{min}$. Then, from equation (9) and (10),

$$\Delta Q = S_i * \sin(\Delta\theta_{min}) \quad (21)$$
$$L_\theta = \sqrt{2}S_i\sqrt{1-\cos(\Delta\theta_n)} \quad (22)$$

### C. Linear Optimal Power Flow Problem

The Linear OPF (LOPF) is formulated here is based on the Decoupled Linearized Power Flow (DLPF) presented in [4], for the generation cost minimization objective. All constraints are linear except the branch flow one. The branch flow constraint:

$$-S_{ij}^2 \leq (P_{ij}^2 + Q_{ij}^2) \leq S_{ij}^2 \quad (23)$$

Where, $i, j$ are bus indices. The number of linear limits for one branch is equal to the number of segments $M$ and given as a set:

$$L_b = \{L_{12}^b, \dots\dots\dots\dots L_{M_b-1,M_b}^b\} \quad (24)$$
$$L_{n,n+1}^b: a_{n,n+1}^b P + b_{n,n+1}^b Q + c_{n,n+1}^b = 0 \quad (25)$$

As $P = f(V,\delta)$ and $Q = h(V,\delta)$ then:
$$L_{n,n+1}^b: a_{n,n+1}^b f(V,\delta) + b_{n,n+1}^b h(V,\delta) + c_{n,n+1}^b \geq 0 \quad (26)$$

$$L = \{L_b\} \geq 0 \quad (27)$$

Here, $b$ is branch index; $n = 1,2 \dots \dots M_b(b)$ and $M_b$ is a vector consisting $M_{irr}$ value for all $b$ branches. Therefore, each equation of branch flow limit (22) will get replaced by a set of equation $L \geq 0$ (27). The number of linear constraints for one branch is equal to the number of sides of the polygon. This is because to retrace the circle it is replaced by the irregular polygon. To do so, $M_{irr}$ line equations, corresponding to the sides of polygon, must be given as linear constraints.

## II. RESULTS

Table I shows the values of $M_{reg}$ and $M_{irr}$ for different errors validating the MVA limit based linearization requirement. The angle inscribed by first line segment for different branch limits is given in Fig. (5). This reflects upon the fact that for same value of $\Delta\theta_1$, lower error can be selected. As the relations (6), (17) and (20) are nonlinear, post rounding off same number of segments can provide lower error as well. The difference in number of branch constraints for same standard system is large enough to provide significant time reduction during solution as shown in Table II.

For the 39-Bus system, LOPF with $e_{m,min}$=0.01 MVA, it is observed that many branches operate in the area which is suggested to be omitted by taking $\alpha$=30º [2], [3] as shown in Fig. (6). This clearly implies that selection of $\alpha$ cannot be done arbitrarily. It reflects that approximation error has a minimum impact because heavily loaded lines operate near P-axis and near the Q-axis, the flow is not restricted by the limiting condition ($Q_{ij}$ not as high as limit in such cases).

The proposed method only changes the branch flow constraints in the linear OPF model. Thus, it will affect the objective value only if it changes the power flow pattern significantly. Fig. (7) shows the normalized power flows with regular and irregular polygon linearization for IEEE-30 bus system and Fig. (8) shows same for IEEE-118 bus system. It depicts that proposed method does not affect the flow pattern as such and thus will not have an impact on optimal value as well. This also shows that proposed approach does not affect the model accuracy of OPF formulation.

The Table III shows the objective value error in comparison to the MATPOWER solution. The results are obtained with $e_{m,reg} = e_{m,min} = 0.01$ MVA. The values indicate that proposed method is as efficient as the regular polygon one in terms of optimality.

As given in Table I and II, the proposed method provides considerable reduction in the number of constraints and the computation results also show a considerable reduction in time efforts of the solution. The internal time taken for IEEE-30 bus systems with different error values ($e_{m,reg}$ or $e_{m,min}$) is shown in Fig. (9) when a system with Intel i5 processor of 3.30GHz clock speed is used with MATLAB R2017a.

The results explain that as we move towards the higher accuracy, the time saving by proposed method is increased. Further, the percentage of time-saving is given in Fig. (10) which emphasizes on improved computational performance of the proposed method.

## III. CONCLUSION

This letter proposes a piecewise linearization method which utilizes operational observations and can be used for thermal limit linearization. Chances of non-convergence and constraint matrix size are reduced. The proposed approach gives an advantage over the existing ones as it provides significant computational time reduction without affecting the optimal value and model accuracy. The work presents a criterion to get irregular polygon and, more criteria can be developed for specific applications, like a dynamic method where highly loaded branches get linearized with low error. The proposed irregular polygon approach is more suitable, in comparison to existing ones, as it gives more flexibility and adaptability for speed-accuracy trade-off. Also, the relation of number of polygon segments with error and MVA limits allows fast yet accurate LOPF solutions.

TABLE I
$M_{reg}$ AND $M_{irr}$ FOR DIFFERENT BRANCH LIMITS

| $e_{m,reg}$ or $e_{m,min}$ (MVA) | Branch Flow Limits (MVA) | | | | | | | |
|---|---|---|---|---|---|---|---|---|
| | Regular Polygon | | | | Irregular Polygon | | | |
| | 16 | 220 | 880 | 1800 | 16 | 220 | 880 | 1800 |
| 0.1 | 29 | 105 | 209 | 299 | 20 | 68 | 136 | 192 |
| 0.2 | 20 | 74 | 148 | 211 | 16 | 48 | 96 | 136 |
| 0.3 | 17 | 61 | 121 | 173 | 12 | 40 | 80 | 112 |

TABLE II
NUMBER OF TOTAL BRANCH CONSTRAINTS FOR DIFFERENT SYSTEMS

| $e_{m,reg}$ or $e_{m,min}$ (MVA) | Standard IEEE Test Systems | | | | | |
|---|---|---|---|---|---|---|
| | Regular Polygon | | | Irregular Polygon | | |
| | 118-Bus | 39-Bus | 30-Bus | 118-Bus | 39-Bus | 30-Bus |
| 0.1 | 20982 | 8684 | 1884 | 13588 | 5620 | 1288 |
| 0.2 | 14805 | 6148 | 1340 | 9608 | 3984 | 952 |
| 0.3 | 12188 | 5018 | 1092 | 7980 | 3260 | 768 |

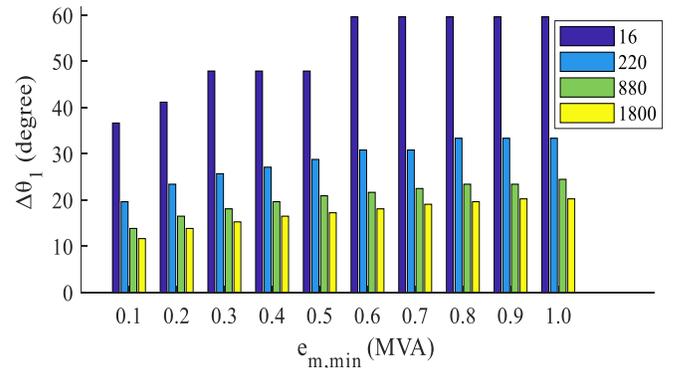

Fig. 5. Variation in $\Delta\theta_1$ with the change in $e_{m,min}$ for different branch limits

TABLE III
ERROR IN OBJECTIVE VALUE AND NUMBER OF INEQUALITY CONSTRAINTS

| System | Error in Objective Value ($/hr.) | | Number of Inequality Constraints | |
|---|---|---|---|---|
| | Irregular Polygon | Regular Polygon | Irregular Polygon | Regular Polygon |
| 9-Bus | 0.36% | 0.57% | 672 | 1018 |
| 30-Bus | 0.70% | 0.72% | 3348 | 5146 |
| 39-Bus | 2.53% | 2.53% | 5798 | 8898 |
| 118-Bus | 1.90% | 1.92% | 36700 | 56728 |

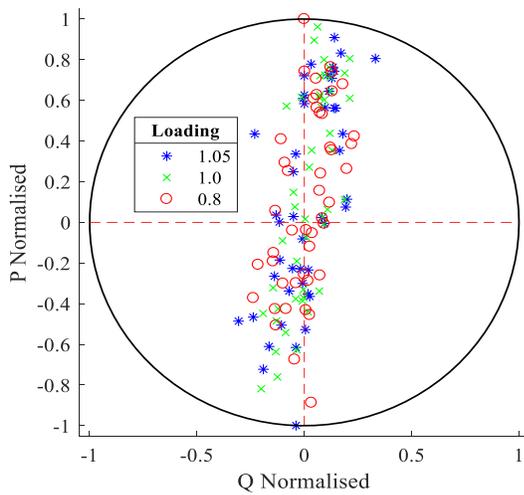

Fig. 6. Normalized value of line flows for 39-Bus system in P-Q plane.

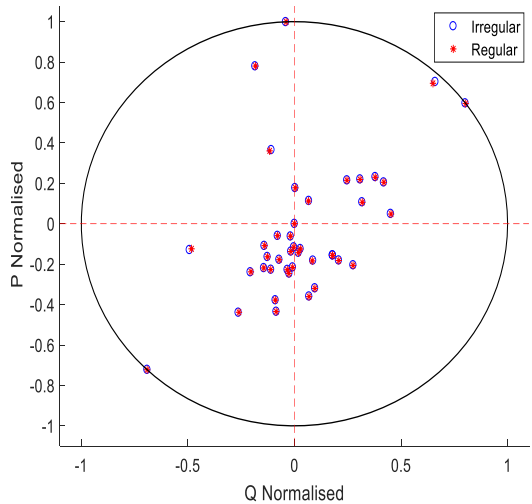

Fig. 7. Normalized power flow for IEEE-30 bus system

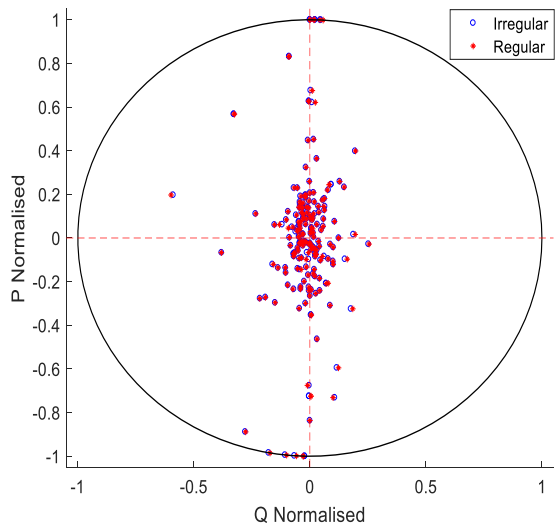

Fig. 8. Normalized power flow for IEEE-118 bus system

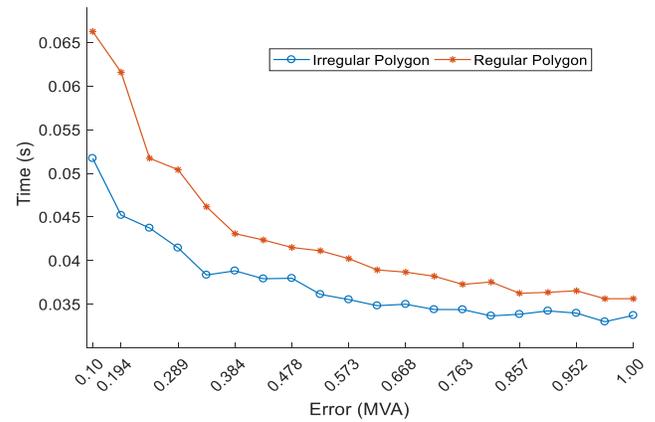

Fig. 9. Internal time taken for OPF solution for IEEE-30 bus system for different values of error ($e_{m,reg}$ or $e_{m,min}$)

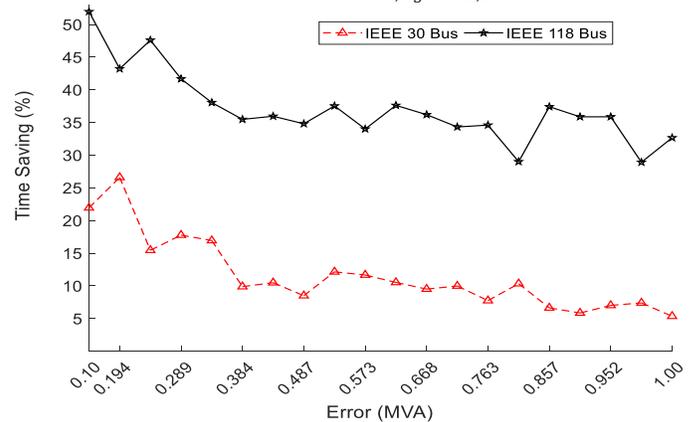

Fig. 10. Percentage time savings by proposed method for different values of error ($e_{m,reg}$ or $e_{m,min}$)

## References

0

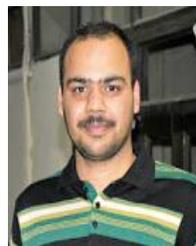

**Parikshit Pareek** was born in India in 1993. He has completed the M.Tech. in Energy Studies from Indian Institute of Technology Delhi, India in 2018. He is selected for Shirimati Jawala Devi-Sita Ram Kaushik Award for sesstion 2017-18 for best M.Tech. project in the Centre for Energy Studies.
He recently joined School of Electrical and Electronics Engineering, Nanyang Technological University, Singapore for Ph.D. His area of interest is power system optimization and applications.

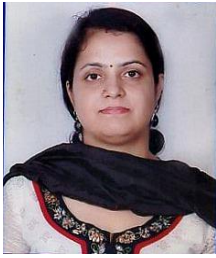

**Ashu Verma** (M'12, SM'17) received the M.Tech. degree in power systems and the Ph.D. degree in transmission expansion planning from the Indian Institute of Technology Delhi (IIT Delhi), New Delhi, India, in 2002 and 2010, respectively.

She is an Assistant Professor with the Centre for Energy Studies, IIT Delhi. She is heading the Electrical Power and Renewable Energy Systems (EPRES) Group at CES, IIT Delhi. Her present research interests include power system planning, operation and control aspects of integrated renewable energy systems, and policy and regulatory framework for enabling large-scale integration of renewable energy sources in India.